\documentclass[12pt]{article}
\usepackage[latin1]{inputenc}
\usepackage{amsmath}
\usepackage{amsfonts}
\usepackage{amssymb}
\usepackage{amsthm}

\usepackage{fixltx2e} 

\usepackage{titlesec}
\titleformat{\section}[hang]{\normalfont\large\bfseries}{\thesection}{1em}{}
\titleformat{\subsection}[runin]{\normalfont\bfseries}{\thesubsection}{1.5ex}{}

\titlespacing\section{0pt}{3.5ex plus 0.5ex minus .2ex}{0.3ex plus .2ex}
\titlespacing\subsection{0pt}{2.5ex plus 0.5ex minus .2ex}{1em}
\titlespacing\subsubsection{0pt}{2.5ex plus 0.5ex minus .2ex}{0.3ex plus .2ex}

\usepackage{mathrsfs} 

\usepackage[alphabetic,lite]{amsrefs} 

\usepackage{fullpage}

\usepackage{setspace}
\usepackage{hyperref}
\usepackage{color}
\usepackage{enumitem}
\usepackage{ulem} 

\usepackage{comment} 

\usepackage[all,cmtip]{xy} 

 \usepackage{fancyhdr}
\newtheorem{Thm}[subsection]{Theorem}

\newtheorem{Lemma}[subsection]{Lemma}

\newtheorem{Cor}[subsection]{Corollary}

\newtheorem{Prop}[subsection]{Proposition}

\theoremstyle{definition}

\newtheorem{Rem}[subsection]{Remark}

\newcommand{\Proof}{\textbf{Proof.\\}}

\usepackage{marginnote}

\newcommand{\<}{\left\langle}
\renewcommand{\>}{\right\rangle}

\newcommand{\bC}{\mathbb{C}}

\newcommand{\bF}{\mathbb{F}}
\newcommand{\bG}{\mathbb{G}}

\newcommand{\bR}{\mathbb{R}}

\newcommand{\bZ}{\mathbb{Z}}

\newcommand{\cO}{\mathcal{O}}
\newcommand{\cP}{\mathcal{P}}

\newcommand{\fg}{\mathfrak{g}}

\newcommand{\fr}{\mathfrak{r}}

\newcommand{\ft}{\mathfrak{t}}

\newcommand{\sA}{\mathscr{A}}
\newcommand{\sB}{\mathscr{B}}

\DeclareMathAlphabet{\mathpzc}{OT1}{pzc}{m}{it}

\newcommand{\ra}{\rightarrow}

\newcommand{\wt}{\widetilde}

\newcommand{\ov}{\overline}

\DeclareMathOperator{\Hom}{Hom}			
\DeclareMathOperator{\Id}{Id}		  	
\DeclareMathOperator{\GL}{GL}		  	
\DeclareMathOperator{\U}{U} 			
\DeclareMathOperator{\Ind}{Ind}			
\DeclareMathOperator{\Sp}{Sp}		  	
\DeclareMathOperator{\diag}{diag}  	
\DeclareMathOperator{\Lie}{Lie}  	  
\DeclareMathOperator{\ind}{ind}  	  
\DeclareMathOperator{\cind}{c-ind}  	  
\DeclareMathOperator{\val}{val}  	  
\DeclareMathOperator{\Cent}{Cent}  	  
\DeclareMathOperator{\pr}{pr} 		

\newcommand{\ba}{\begin{aligned}}
\newcommand{\ea}{\end{aligned}}
\newcommand{\beqn}{\begin{eqnarray}}
\newcommand{\eeqn}{\end{eqnarray}}
\newcommand{\beqns}{\begin{eqnarray*}}
\newcommand{\eeqns}{\end{eqnarray*}}
\newcommand{\benum}{\begin{enumerate}}
\newcommand{\eenum}{\end{enumerate}}

\newcommand{\repKYu}{\wt \rho}

\newcommand{\KYu}{\wt K}  	  
\newcommand{\Uff}{\bar U}  	 
\newcommand{\UU}{U} 
\newcommand{\Uiff}{\bar U_{i}}

\newcommand{\Htt}{H_{23}} 
\newcommand{\Hff}{H_{45}} 
\newcommand{\fsp}{\mathfrak{sp}} 
 
\newcommand{\tildephi}{\widetilde \phi}

\AtEndDocument{\bigskip{\footnotesize%
		\textsc{University of Cambridge, Cambridge, UK and Duke University, Durham, NC, USA} \par  
		\textit{Mailing address}: Trinity College, Cambridge, CB2 1TQ, UK \par
		\textit{E-mail address}: \texttt{fintzen@maths.cam.ac.uk} and \texttt{fintzen@math.duke.edu}
}}

\begin{document}
\author{Jessica Fintzen}
\title{On the construction of tame supercuspidal representations} 
\date{}
\maketitle
\begin{abstract}
  Let $F$ be a non-archimedean local field of residual characteristic $p \neq 2$. Let $G$ be a (connected) reductive group over $F$ that splits over a tamely ramified field extension of $F$. We revisit Yu's construction of smooth complex representations of $G(F)$ from a slightly different perspective
   and provide a proof that the resulting representations are supercuspidal. 
   
   We also provide a counterexample to Proposition 14.1 and Theorem 14.2 in \cite{Yu}, whose proofs relied on a typo in a reference.
   
  \end{abstract}

{
	\renewcommand{\thefootnote}{}  
	\footnotetext{MSC2010: 22E50; 20G05, 20G25} 
	\footnotetext{Keywords: representations of reductive groups over non-archimedean local fields, supercuspidal representations, $p$-adic groups}
	\footnotetext{The author was partially supported by NSF Grant DMS-1802234 / DMS-2055230 and a Royal Society University Research Fellowship.}
}

\setcounter{tocdepth}{1}

\tableofcontents

\newpage

\section{Introduction}
In 2001, Yu (\cite{Yu}) proposed a construction of smooth complex supercuspidal representations of $p$-adic groups that since then has been widely used, e.g. to study the Howe correspondence, to understand distinction of representations of $p$-adic groups, to obtain character formulas and to construct an explicit local Langlands correspondence. However, Loren Spice noticed recently that Yu's proof relies on a misprinted\footnote{As Loren Spice pointed out, the statement of \cite[Theorem~2.4.(b)]{Gerardin} contains a typo. From the proof provided by \cite{Gerardin} one can deduce that the stated representation of $P(E_+,j)H(E_+^\perp,j)$ (i.e. the pull-back to $P(E_+,j)H(E_+^\perp,j)$ of a representation of $SH(E_0,j_0)$ as in part (a')) should be tensored with $\chi^{E_+} \ltimes 1$ before inducing it to $P(E_+,j)H(E,j)$ in order to define $\pi_+$ (using the notation of \cite{Gerardin}).} (and therefore false) statement in \cite{Gerardin} and therefore it became uncertain whether the representations constructed by Yu are irreducible and supercuspidal.
In the present paper we illustrate the significance of this false statement on Yu's proof by providing a counterexample to Proposition 14.1 and Theorem 14.2 of \cite{Yu}. Proposition 14.1 and Theorem 14.2 are the main intertwining results in \cite{Yu} that form the heart of the proof. We also offer a different argument to show that nevertheless Yu's construction yields irreducible supercuspidal representations. 

 Let $F$ be a non-archimedean local field of residual characteristic $p \neq 2$. Let $G$ be a (connected) reductive group that splits over a tamely ramified field extension of $F$. In this paper, we first describe the construction of Yu's representations in a way that we find more convenient for our purpose and then provide a proof that these representations are supercuspidal. All representations arise via compact induction from an irreducible representation $\repKYu$ of a compact-mod-center open subgroup $\KYu$ of $G(F)$. 
 Our proof only relies on the first part of Yu's proof and provides a shorter, alternative second part that does not rely on \cite[Proposition~14.1~and~Theorem~14.2]{Yu} and the misprinted version of \cite[Theorem~2.4(b)]{Gerardin}. Yu's approach consists of following a strategy already employed by Bushnell--Kutzko that required to show that a certain space of intertwining operators has dimension precisely one, i.e., in particular, is non-trivial. Our approach does not require such a result. Instead we use the structure of the constructed representation including the structure of Weil--Heisenberg representations, and the Bruhat--Tits building to show more directly that every element that intertwines $\repKYu$ is contained in $\KYu$, which implies the desired result. Our proof  
  relies also less heavily on tameness assumptions, and our aim is to use a modification of it for the construction of supercuspidal representations beyond the tame setting when Yu's construction is not exhaustive.

 Note that Yu's construction yields all supercuspidal representations if $p$ does not divide the order of the Weyl group of $G$ (\cites{Fi-exhaustion, Kim}), a condition that guarantees that all tori of $G$ split over a tamely ramified field extension of $F$.
 
 In the last section we provide a counterexample to \cite[Proposition~14.1~and~Theorem~14.2]{Yu} by considering the group $G=\Sp_{10}$ together with a twisted Levi subgroup $G'$ of shape $\U(1) \times \Sp_8$ and a well chosen point in the Bruhat--Tits building of $G'$. 

\textbf{Conventions and notation.} 
Let $F$ be a non-archimedean local field of residual characteristic $p \neq 2$. We denote by $\cO$ the ring of integers of $F$, and by $\cP$ the maximal ideal of $\cO$. The residue field $\cO/\cP$ is denoted by $\bF_q$, where $q$ denotes the number of elements in $\bF_q$. When considering field extensions of $F$ in this paper, we mean algebraic field extensions of $F$ and view them as contained in a fixed algebraic closure $\ov F$ of $F$. If $E$ is a field extension of $F$, then we write $E^{ur}$ for the maximal unramified extension of $E$ and $E^{sep}$ for the separable closure of $E$. 

All reductive groups in this paper are required to be connected.

For a reductive group $G$ defined over $F$ we denote by $\sB(G,F)$ the (enlarged) Bruhat--Tits building (\cites{BT1,BT2}) of $G$ over $F$, by $Z(G)$ the center of $G$ and by $G^{\mathrm{der}}$ the derived subgroup of $G$. If $T$ is a maximal, maximally split torus of $G_E:=G \times_F E$ for some field extension $E$ over $F$, then $\sA(T,E)$ denotes the apartment of $T$ inside the Bruhat--Tits building $\sB(G_E,E)$ of $G_E$ over $E$. Moreover, we write $\Phi(G_E,T)$ for the roots of $G_{E} \times_{E} \ov F$ with respect to $T_{\ov F}$.  We let $\wt \bR=\bR \cup \{ {r+} \, | \, r \in \bR\}$ with its usual order, i.e. for $r$ and $s$ in $\bR$ with $r<s$, we have $r<r+<s<s+$. For $r \in \wt \bR_{\geq 0}$, we write $G_{x,r}$ for the Moy--Prasad filtration subgroup of $G(F)$ of depth $r$ at a point $x \in \sB(G,F)$. For $r \in \wt \bR$, we write $\fg_{x,r}$ for the Moy--Prasad filtration submodule of  $\fg=\Lie G(F)$ of depth $r$ at $x$, and $\fg^*_{x,r}$ for the Moy--Prasad filtration submodule of depth $r$ at $x$ of the linear dual $\fg^*$ of $\fg$.  If $x \in \sB(G,F)$, then we denote by $[x]$ its image in the reduced Bruhat--Tits building and we write 
 $G_{[x]}$ for the stabilizer of $[x]$ in $G(F)$. 

We call a subgroup $G'$ of $G$ (defined over $F$) a twisted Levi subgroup of $G$ if $(G')_E$ is a Levi subgroup of $G_E$ for some (finite) field extension $E$ of $F$. If $G'$ splits over a tamely ramified field extension of $F$, then, using (tame) Galois descent, we obtain an embedding of the corresponding Bruhat--Tits buildings $\sB(G',F) \hookrightarrow \sB(G,F)$. This embedding is only unique up to some translation, but its image is unique, and we will identify $\sB(G',F)$ with its image in $\sB(G,F)$. All constructions in this paper are independent of the choice of such an identification.

Let $H$ be a group and $\chi$ a character of $H$. Then we denote by $\bC_\chi$ the one dimensional complex representation space on which $H$ acts via $\chi$. We also write $1$ to denote the one dimensional trivial complex representation.
If $K$ is a subgroup of $H$, $h \in H$, and $\rho$ a representation of $K$, then we write $^hK$ to denote $hKh^{-1}$ and define ${^h\rho}(x)=\rho(h^{-1}xh)$ for $x \in K \cap {^hK}$. We say that $h$ intertwines $\rho$ if the space of intertwiners $\Hom_{K \cap {^hK}}(\rho, {^h\rho})$ is non-zero.

Throughout the paper we fix an additive character $\varphi: F \ra \bC^*$ of $F$ of conductor $\cP$ and a reductive group $G$ that is defined over our non-archimedean local field $F$ and that splits over a tamely ramified field extension of $F$. All representations of $G(F)$ in this paper have complex coefficients and are required to be smooth.

\textbf{Acknowledgments.} The author thanks Loren Spice for pointing out that Yu's proof relies on a misprinted (and therefore false) statement in a paper by Gérardin, Tasho Kaletha for his encouragement to write up the below presented proof that Yu's construction yields irreducible supercuspidal representations, and the referees for helpful feedback on the paper.

The author also thanks Jeffrey Adler, Stephen DeBacker, Tasho Kaletha, Loren Spice and Cheng-Chiang Tsai for helpful discussions related to the topic of this paper during their SQuaRE meetings at the American Institute of Mathematics, as well as the American Institute of Mathematics for supporting these meetings and providing a wonderful research environment.

\section{Construction of representations à la Yu}\label{Section-input-and-construction}

In this section we recall Yu's construction of representations but formulated in a way that is better adapted to our proof of supercuspidality. 

\subsection{The input} \label{Section-input}

The input for Yu's construction of supercuspidal representations of $G(F)$ (using the conventions from \cite{Fi-exhaustion}, see Remark \ref{Rem-conventions} for a comparison of Yu's notation with ours) 
 is a tuple $((G_i)_{1 \leq i \leq n+1}, x, (r_i)_{1 \leq i \leq n}, \rho, (\phi_i)_{1 \leq i \leq n})$ for some non-negative integer $n$ where
		\begin{enumerate}[label=(\alph*),ref=\alph*]
			\item $G=G_1 \supseteq G_2 \supsetneq G_3 \supsetneq \hdots \supsetneq G_{n+1}$ are twisted Levi subgroups of $G$ that split over a tamely ramified extension of $F$
			\item  $x \in \sB(G_{n+1},F)\subset \sB(G,F)$ 
			\item $r_1 > r_2 > \hdots > r_n >0$ are real numbers
			\item \label{item-rho} $\rho$ is an irreducible representation of $(G_{n+1})_{[x]}$ that is trivial on $(G_{n+1})_{x,0+}$ 
			\item $\phi_i$, for $1 \leq i \leq n$, is a character of $G_{i+1}(F)$ of depth $r_i$ that is trivial on $(G_{i+1})_{x,r_i+}$ 
		\end{enumerate}
		satisfying the following conditions 
		\begin{enumerate}[label=(\roman*),ref=\roman*]
			\item \label{item-anisotropic} $Z(G_{n+1})/Z(G)$ is anisotropic
			\item \label{condition-vertex} the image of the point $x$ in $\sB(G_{n+1}^{\mathrm{der}},F)$ is a vertex
			\item \label{condition-cuspidal} $\rho|_{(G_{n+1})_{x,0}}$ is a cuspidal representation of $(G_{n+1})_{x,0}/(G_{n+1})_{x,0+}$
			\item \label{condition-generic} $\phi_i$ is $G_i$-generic of depth $r_i$ relative to $x$ (in the sense of \cite[\S9, p.~599]{Yu}) for all $1 \leq i \leq n$ with $G_i \neq G_{i+1}$ .
		\end{enumerate}

\begin{Rem} \label{Rem-tame-torus}
	Note that for each apartment $\sA$ of $\sB(G_{n+1},F)$, there exists a maximal torus $T$ of $G_{n+1}$ that splits over a tamely ramified extension $E$ of $F$ such that $\sA \subset \sA(T, E)$ (see, e.g., \cite[\S2, page~585-586]{Yu}, which is based on \cite{BT2} and \cite{Rousseau}). In particular, there exists a maximal torus $T$ of $G_{n+1}$ that splits over a tamely ramified extension $E$ of $F$ such that $x \in \sA(T, E)$.
\end{Rem}

\begin{Rem} \label{Rem-rho}
	By (the proof of) \cite[Proposition~6.8]{MP2} requiring that the image of the point $x$ in $\sB(G_{n+1}^{\mathrm{der}},F)$ is a vertex and that $\rho|_{(G_{n+1})_{x,0}}$ is a cuspidal representation of $(G_{n+1})_{x,0}/(G_{n+1})_{x,0+}$ is equivalent to requiring that $\cind_{(G_{n+1})_{[x]}}^{G_{n+1}(F)}\rho$ is an irreducible supercuspidal representation.
	When $n=0$, then the tuple $((G_i)_{1 \leq i \leq n+1}, x, (r_i)_{1 \leq i \leq n}, \rho, (\phi_i)_{1 \leq i \leq n})$ consists only of the group $G=G_1=G_{n+1}$, a point $x \in \sB(G,F)$ whose image in $\sB(G^{\mathrm{der}},F)$ is a vertex and an irreducible representation $\rho$ of $G_{[x]}$ that is trivial on $G_{x,0+}$ and such that its restriction $\rho|_{G_{x,0}}$ is a cuspidal representation of $G_{x,0}/G_{x,0+}$. This case recovers the depth-zero supercuspidal representations.
\end{Rem}

\begin{Rem} \label{Rem-conventions}
	We use the conventions for notation from \cite{Fi-exhaustion} instead of from \cite{Yu}. The notation in \cite{Yu} (left hand side) can be recovered from ours (right hand side) as follows:	
	\begin{eqnarray*} 
		\vec G= (G^0, G^1, \hdots, G^d) &=& \left\{ \begin{array}{ll} (G_{n+1}, G_n, \hdots, G_2, G_1=G) & \text{ if } G_2 \neq G_1 \text{ or } n=0 \\
			(G_{n+1}, G_n, \hdots, G_3, G_2=G) & \text{ if } G_2=G_1 \end{array}\right. 
		\\
		\vec r & = & \left\{ \begin{array}{ll} (r_{n}, r_{n-1}, \hdots, r_2, r_1, r_\pi) & \text{ if } G_2 \neq G_1 \text{ or } n=0 \\
			(r_{n}, r_{n-1}, \hdots, r_2, r_1)  & \text{ if } G_2=G_1 \end{array}\right.
		\\
		\vec \phi & = & \left\{ \begin{array}{ll} (\phi_{n}, \phi_{n-1}, \hdots, \phi_2, \phi_1, 1) & \text{ if } G_2 \neq G_1 \text{ or } n=0 \\
			(\phi_{n}, \phi_{n-1}, \hdots, \phi_2, \phi_1)  & \text{ if } G_2=G_1 \end{array}\right. 	,
	\end{eqnarray*}
where $r_\pi=r_1$ if $n \geq 1$ and $r_\pi=0$ if $n=0$.
Yu's convention has the advantage that it is adapted to associating a whole sequence of supercuspidal representations to a given datum (by only considering the groups $G_i, G_{i+1}, \hdots, G_{n+1}$), while our convention is more natural when recovering the input from a given representation as can be seen in \cite{Fi-exhaustion}. We have chosen our convention for this paper as it has the advantage that our induction steps below start with $G_1$ and move from $G_i$ to $G_{i+1}$. Moreover, using our notation we do not have to impose a condition on $\phi_d$ depending on whether $r_{d-1} < r_d$ or $r_{d-1}=r_d$ in Yu's notation, see \cite[page~590~D5]{Yu}. Hence the input looks more uniform.
(Note that our condition $G_i \neq G_{i+1}$ in \eqref{condition-generic} could be removed by extending the notion of $G_i$-generic to the case $G_i = G_{i+1}$.)
\end{Rem}

\subsection{The construction}\label{Section-construction}

The (smooth complex) representation $\pi$ of $G(F)$ that Yu constructs from the given input  $((G_i)_{1 \leq i \leq n+1}, x, (r_i)_{1 \leq i \leq n}, \rho, (\phi_i)_{1 \leq i \leq n})$ is the compact induction $\cind_{\KYu}^{G(F)} \repKYu$ of a representation $\repKYu$ of a compact-mod-center, open subgroup $\KYu \subset G(F)$. 

In order to define $\KYu$ and $\repKYu$ we introduce the following notation. For $\wt r \geq \wt r' \geq \frac{\wt r}{2}>0$ ($\wt r, \wt r' \in \wt \bR$) and $1 \leq i \leq n$, we choose a maximal torus $T$ of $G_{i+1}$ that splits over a tamely ramified extension $E$ of $F$ and such that $x \in \sA(T,E)$. Then we define
\begin{eqnarray*}
	 (G_{i})_{x,\wt r,\wt r'} &:=& G(F) \cap \\
	  && \<T(E)_{\wt r}, U_\alpha(E)_{x,\wt r}, U_\beta(E)_{x,\wt r'} \, | \, \alpha \in \Phi(G_i, T)\subset\Phi(G,T), \beta \in \Phi(G_i, T)-\Phi(G_{i+1}, T) \, \>,
\end{eqnarray*}
where 
$U_\alpha(E)_{x,r}$ denotes the Moy--Prasad filtration subgroup of depth $r$ (at $x$) of the root group $U_\alpha(E) \subset G(E)$ corresponding to the root $\alpha$. We define $(\fg_{i})_{x,\wt r,\wt r'}$ analogously for $\fg_i=\Lie(G_i)(F)$. The group $(G_{i})_{x,\wt r, \wt r'}$ is denoted by $(G_{i+1},G_i)(F)_{x_i,\wt r, \wt r'}$ in \cite{Yu}, and Yu (\cite[p.~585 and p.~586]{Yu}) shows that this definition is independent of the choice of $T$ and $E$.

We set 
\begin{eqnarray*}
	\KYu &=&(G_1)_{x,\frac{r_1}{2}}(G_2)_{x,\frac{r_2}{2}} \hdots (G_n)_{x,\frac{r_n}{2}}(G_{n+1})_{[x]} \\
	&=&(G_1)_{x,r_1,\frac{r_1}{2}}(G_2)_{x,r_2,\frac{r_2}{2}} \hdots (G_n)_{x,r_n,\frac{r_n}{2}}(G_{n+1})_{[x]}. 
\end{eqnarray*}
Note that since we assume that $Z(G_{n+1})/Z(G)$ is anisotropic (see Condition (\ref{item-anisotropic})), the subgroup $\KYu$ of $G(F)$ is compact mod center.
Now the representation $\repKYu$ of $\KYu$ is given by $\rho \otimes \kappa$, where $\rho$ also denotes the extension of $\rho$ from $(G_{n+1})_{[x]}$ to $\KYu$ that is trivial on $(G_1)_{x,\frac{r_1}{2}}(G_2)_{x,\frac{r_2}{2}} \hdots (G_n)_{x,\frac{r_n}{2}}$. In order to define $\kappa$ we need some additional notation.

Following \cite[\S~4]{Yu}, we denote by $\hat \phi_i$ \label{page-hatphi} for $1 \leq i \leq n$ the unique character of $(G_{n+1})_{[x]}(G_{i+1})_{x,0}G_{x,\frac{r_i}{2}+}$ that satisfies
\begin{itemize}
	\item $\hat \phi_i|_{(G_{n+1})_{[x]}(G_{i+1})_{x,0}}=\phi_i|_{(G_{n+1})_{[x]}(G_{i+1})_{x,0}}$, and 
	\item $\hat \phi_i|_{G_{x,\frac{r_i}{2}+}} $ factors through 
	\begin{eqnarray*} 
		G_{x,\frac{r_i}{2}+}/G_{x,r_i+} & \simeq& \fg_{x,\frac{r_i}{2}+}/\fg_{x,r_i+} = (\fg_{i+1} \oplus \fr'')_{x,\frac{r_i}{2}+}/(\fg_{i+1} \oplus \fr'')_{x,r_i+} \\
		& \ra & (\fg_{i+1})_{x,\frac{r_i}{2}+}/(\fg_{i+1})_{x,r_i+} \simeq 
		(G_{i+1})_{x,\frac{r_i}{2}+}/(G_{i+1})_{x,r_i+},
	\end{eqnarray*}
	on which it is induced by $\phi_i$. Here $\fr''$ is defined to be $\fg \cap \bigoplus_{\alpha \in \Phi(G,T_E)-\Phi(G_{i+1},T_E)} (\fg_E)_\alpha$ for some maximal torus $T$ of $G_{i+1}$ that splits over a tame extension $E$ of $F$ with $x \in \sA(T,E)$, and the surjection $\fg_{i+1} \oplus \fr'' \twoheadrightarrow \fg_{i+1}$ sends $\fr''$ to zero. (Recall that $\fg=\Lie(G)(F)$, and $(\fg_E)_{\alpha}$ denotes the $E$-subspace of $\Lie(G)(E)$ on which the torus acts via $\alpha$.)
\end{itemize}

Note that $(G_i)_{x,r_i,\frac{r_i}{2}}/\left((G_i)_{x,r_i,\frac{r_i}{2}+} \cap \ker\hat \phi_i\right)$ is a Heisenberg $p$-group with center \linebreak $(G_i)_{x,r_i,\frac{r_i}{2}+} / \left((G_i)_{x,r_i,\frac{r_i}{2}+} \cap \ker\hat \phi_i\right)$ (\cite[Proposition~11.4]{Yu}).
More precisely, set $$V_i:=(G_i)_{x,r_i,\frac{r_i}{2}}/(G_i)_{x,r_i,\frac{r_i}{2}+} ,$$
and equip it with the pairing $\<\cdot, \cdot\>_i$ defined by $\<a,b\>_i=\hat \phi_i(aba^{-1}b^{-1})$. Then Yu shows in (\cite[Proposition~11.4]{Yu}) that there is a canonical  special isomorphism $$j_i:(G_i)_{x,r_i,\frac{r_i}{2}}/\left((G_i)_{x,r_i,\frac{r_i}{2}+} \cap \ker\hat \phi_i\right) \ra V_i^\sharp ,$$
 where $V_i^\sharp$ is the group with underlying set $V_i \times \bF_p$ and with group law $(v,a).(v',a')=(v+v', a+a'+\frac{1}{2}\<v,v'\>_i)$.
  
Let $(\omega_i, V_{\omega_i})$ denote the Heisenberg representation of $(G_i)_{x,r_i,\frac{r_i}{2}}/\left((G_i)_{x,r_i,\frac{r_i}{2}+} \cap \ker\hat \phi_i\right)$ (via the above special isomorphism) with central character $\hat \phi_i|_{(G_i)_{x,r_i,\frac{r_i}{2}+} }$. Then we define the space $V_\kappa$ underlying the representation $\kappa$ to be $\bigotimes_{i=1}^n V_{\omega_i}$. If $n=0$, then the empty tensor product should be taken to be a one dimensional complex vector space and $\kappa$ is the trivial representation. 
In order to describe the action of $\KYu$ on each $V_{\omega_i}$ for $n \geq 1$, we describe the action of $(G_i)_{x,r_i,\frac{r_i}{2}}$ for $1  \leq i \leq n$ and of $(G_{n+1})_{[x]}$ separately.

For $1  \leq i \leq n$, the action of $(G_i)_{x,r_i,\frac{r_i}{2}}$ on $V_{\omega_i}$ should be given by letting $(G_i)_{x,r_i,\frac{r_i}{2}}$ act via the Heisenberg representation $\omega_i$ of $(G_i)_{x,r_i,\frac{r_i}{2}}/\left((G_i)_{x,r_i,\frac{r_i}{2}+} \cap \ker\hat \phi_i\right)$ with central character $\hat \phi_i|_{(G_i)_{x,r_i,\frac{r_i}{2}+} }$. The action of $(G_i)_{x,r_i,\frac{r_i}{2}}$  on $V_{\omega_j}$ for $j \neq i$ should be via the character $\hat \phi_j|_{(G_i)_{x,r_i,\frac{r_i}{2}}}$ (times identity).

The action of $(G_{n+1})_{[x]}$ on $V_{\omega_i}$ for $1 \leq i \leq n$ is given by $\phi_i|_{(G_{n+1})_{[x]}}$ times the following representation:  \label{page-tildephi} Let $(G_{n+1})_{[x]}/(G_{n+1})_{x,0+}$ act on $V_{\omega_i}$ by mapping $(G_{n+1})_{[x]}/(G_{n+1})_{x,0+}$ to the symplectic group  $\Sp(V_i)$ of the corresponding symplectic $\bF_p$-vector space $V_i=(G_i)_{x,r_i,\frac{r_i}{2}}/(G_i)_{x,r_i,\frac{r_i}{2}+}$ with pairing 
$\<a,b\>_i=\hat \phi_i(aba^{-1}b^{-1})$ (after choosing a from now on fixed isomorphism between the $p$-th roots of unity in $\bC^*$ and $\bF_p$) and composing this map with the Weil representation (defined in \cite{Gerardin}). 
Here the map from $(G_{n+1})_{[x]}/(G_{n+1})_{x,0+}$ to $\Sp(V_i)$ is induced by the conjugation action of $(G_{n+1})_{[x]}$ on $(G_i)_{x,r_i,\frac{r_i}{2}}$, which (together with the special isomorphism $j_i$) yields a symplectic action in the sense of \cite[\S10]{Yu} by \cite[Proposition~11.4]{Yu}. 

Then the resulting actions of $(G_i)_{x,r_i,\frac{r_i}{2}}$ for $1  \leq i \leq n$ and $(G_{n+1})_{[x]}$ agree on the intersections and hence yield a representation $\kappa$ of $\KYu$ on the space $V_\kappa$.

The representation $\pi=\cind_{\KYu}^{G(F)} \rho \otimes \kappa$ is the smooth representation of $G(F)$ that Yu attaches to the tuple $((G_i)_{1 \leq i \leq n+1}, x, (r_i)_{1 \leq i \leq n}, \rho, (\phi_i)_{1 \leq i \leq n})$, and we prove in the next section that $\pi$ is an irreducible, supercuspidal representation.

\section{Proof that the representations are supercuspidal} \label{Section-supercuspidal}

We keep the notation from the previous section to prove the following theorem in this section.

\begin{Thm} \label{Thm-main}
	The representation $\cind_{\KYu}^{G(F)} \repKYu$ is irreducible, hence supercuspidal.
\end{Thm}

\begin{Rem} \label{Rem-Yu-mistake}
	This theorem follows from \cite[Theorem~15.1]{Yu}. However, the proof in \cite{Yu} relies on \cite[Theorem~2.4(b)]{Gerardin} and  unfortunately the statement of \cite[Theorem~2.4(b)]{Gerardin} contains a typo as Loren Spice pointed out. Therefore Proposition~14.1 and Theorem~14.2 of \cite{Yu}, on which Yu's proof relies, are no longer true. We provide a counterexample in Section \ref{Section-counterexample}.
	
	 Here we use an alternative and shorter approach to prove Theorem \ref{Thm-main} that uses ideas from the first part of Yu's paper (\cite[Theorem~9.4]{Yu}), but that avoids the second part that relies on the misprinted version of the theorem in \cite{Gerardin}. In particular, we do not use \cite[Proposition~14.1 and Theorem~14.2]{Yu}. 
\end{Rem}

In order to show that $\cind_{\KYu}^{G(F)} \repKYu$ is irreducible, we first observe that $\repKYu$ is irreducible.

\begin{Lemma} \label{Lemma-irreducible}
	The representation $\repKYu$ of $\KYu$ is irreducible.
\end{Lemma}
\Proof For $1 \leq i \leq n$ set $K_i = (G_1)_{x,r_1,\frac{r_1}{2}}(G_2)_{x,r_2,\frac{r_2}{2}} \hdots (G_i)_{x,r_i,\frac{r_i}{2}}$ and $K_0=\{1\}$.
We first prove by induction on $i$ that $\otimes_{j=1}^i V_{\omega_j}$ is an irreducible representation of $K_i$ via the action described in Section \ref{Section-construction}. For $i=0$, we take  $\otimes_{j=1}^i V_{\omega_i}$  to be the trivial one dimensional representation and the statement holds. Now assume the induction hypothesis that  $\otimes_{j=1}^{i-1} V_{\omega_{j}}$ is an irreducible representation of $K_{i-1}$.  Suppose $V' \subset \left(\otimes_{j=1}^{i-1} V_{\omega_{j}}\right) \otimes V_{\omega_i}$ is a non-trivial subspace that is $K_i$-stable. Since $K_{i-1}$ acts on $V_{\omega_i}$ via a character (times identity), the subspace $V'$ has to be of the form $\left(\otimes_{j=1}^{i-1} V_{\omega_{j}}\right)  \otimes V'' $ for a $K_i$-stable non-trivial subspace $V''$ of $V_{\omega_i}$. However, since Heisenberg representations are irreducible, $V_{\omega_i}$ is irreducible as a representation of $(G_i)_{x,r_i,\frac{r_i}{2}} \subset K_i$, and therefore $V''=V_{\omega_i}$. Thus $\otimes_{j=1}^i V_{\omega_j}$ is an irreducible representation of $K_i$, and by induction the representation $\kappa$  is an irreducible representation of $K_n$.

Since $K_n$ acts trivially on $\rho$, every irreducible $\KYu$-subrepresentation of $\repKYu=\rho \otimes \kappa$ has to be of the form $\rho' \otimes \kappa$ for an irreducible subrepresentation $\rho'$ of $\rho$. As $\rho$ is irreducible when restricted to $(G_{n+1})_{[x]} \subset \KYu$, we deduce that $\repKYu$ is an irreducible representation of $\KYu$.
\qed

The remaining proof of Theorem \ref{Thm-main} is concerned with showing that if $g$ intertwines $\repKYu$, then $g \in \KYu$, which then implies that $\ind_{\KYu}^{G(F)}\repKYu$ is irreducible and hence supercuspidal. Our proof consists of two parts. The first part is concerned with reducing the problem to considering $g \in G_{n+1}(F)$ using the characters $\phi_i$ and is essentially \cite[Corollary~4.5]{Yu}. The second part consists of deducing from there the theorem using the depth-zero representation $\rho$ together with the action of suitably chosen subgroups of higher depth and employing knowledge about the structure of Weil--Heisenberg representations. This is where our approach deviates crucially from Yu's approach by avoiding the wrong \cite[Proposition~14.1 and Theorem~14.2]{Yu} in favor of a much shorter argument.
For the first part, we will use the following result of Yu (\cite[Theorem~9.4]{Yu}).
\begin{Lemma}[\cite{Yu}] \label{Lemma-Yu} 
	 Let $1 \leq i \leq n$ and $g \in G_i(F)$.
	 Suppose that $g$ intertwines $\hat \phi_i|_{(G_i)_{x,r_i,\frac{r_i}{2}+}}$. Then $g \in (G_i)_{x,\frac{r_i}{2}}G_{i+1}(F)(G_i)_{x,\frac{r_i}{2}}$. 
\end{Lemma}	
\Proof
This is (part of) \cite[Theorem~9.4]{Yu}. \qed

\textbf{Proof of Theorem \ref{Thm-main}.}\\
Recall that $\repKYu$ is irreducible by Lemma \ref{Lemma-irreducible}. Thus, in order to show that  $\cind_{\KYu}^{G(F)} \repKYu$ is irreducible, hence supercuspidal, we have to show that if $g \in G(F)$ such that 
\begin{equation} \label{eqn-intertwiner2}   \Hom_{\KYu \cap {^g\KYu}}\left( ^g\repKYu|_{\KYu \cap {^g\KYu}},  \repKYu|_{\KYu \cap {^g\KYu}}\right) \neq \{0 \} ,
	\end{equation}
then $g \in \KYu$,  where $^g\KYu$ denotes $g\KYu g^{-1}$ and $^g\repKYu(x)=\repKYu(g^{-1}xg)$.

Fix such a $g \in G(F)$ satisfying $ \Hom_{\KYu \cap {^g\KYu}}\left( ^g\repKYu,  \repKYu\right) \neq \{0 \} $, and define
$$\KYu_i = (G_1)_{x,\frac{r_1}{2}}(G_2)_{x,\frac{r_2}{2}} \hdots (G_i)_{x,\frac{r_i}{2}} \quad \text{ and } \quad \KYu_0 = \{ 1 \} .$$
 We first prove by induction that $g \in  \KYu_n G_{n+1}(F)\KYu_n$ using Lemma \ref{Lemma-Yu} (which is (part of) \cite[Theorem~9.4]{Yu}). This is essentially \cite[Corollary~4.5]{Yu}, but we include a short proof for the convenience of the reader. Let $1 \leq i  \leq n$ and assume the induction hypothesis that $g \in  \KYu_{i-1}G_{i}(F)\KYu_{i-1}$, which is obviously satisfied for $i=1$. We need to show that $g \in  \KYu_{i}G_{i+1}(F)\KYu_{i}$. Since $\KYu_{i-1} \subset \KYu_i \subset \KYu$, we may assume without loss of generality that $g \in G_i(F)$.
Recall that by construction $\rho|_{(G_i)_{x,r_i,\frac{r_i}{2}+}}= \Id$ and $\kappa|_{(G_i)_{x,r_i,\frac{r_i}{2}+}}=\prod_{j=1}^n \hat \phi_j \cdot \Id$. 
 Thus by restriction of the action in Equation \eqref{eqn-intertwiner2} to $(G_i)_{x,r_i,\frac{r_i}{2}+} \cap {^{g}(G_i)_{x,r_i,\frac{r_i}{2}+}}$ we conclude that $g$ intertwines $(\prod_{j=1}^n \hat \phi_j)|_{(G_i)_{x,r_i,\frac{r_i}{2}+}}$.
By the definition of $\hat \phi_j$ in Section \ref{Section-construction}, we have that $\hat \phi_j|_{(G_i)_{x,r_i,\frac{r_i}{2}+}}$ is trivial for $j>i$. Moreover, if $j < i$, then for $y\in (G_i)_{x,r_i,\frac{r_i}{2}+} \cap  {^{g}(G_i)_{x,r_i,\frac{r_i}{2}+}}$ we have
\begin{equation} \label{eqn-g'-rho} {^{g}\hat \phi_j}(y)= \hat \phi_j(g^{-1}yg)=\phi_j(g^{-1}yg)=\phi_j(g^{-1})\phi_j(y)\phi_j(g)=\phi_j(y)=\hat \phi_j(y) .
\end{equation}
Therefore we obtain that $g$ also intertwines $\hat \phi_i|_{(G_i)_{x,r_i,\frac{r_i}{2}+}}$.
By Lemma \ref{Lemma-Yu} (which is (part of) \cite[Theorem~9.4]{Yu}) we conclude that $g\in (G_i)_{x,\frac{r_i}{2}}G_{i+1}(F)(G_i)_{x,\frac{r_i}{2}}$, and hence 
$g \in \KYu_i G_{i+1}(F)\KYu_i $. This finishes the induction step and therefore we have shown that $g \in  \KYu_nG_{n+1}(F)\KYu_n$.

In order to prove that $g \in \KYu$, we may therefore assume without loss of generality that $g \in G_{n+1}(F)$, and it suffices to prove that then $g \in (G_{n+1})_{[x]}$. Let us assume the contrary, i.e. $g \in G_{n+1}(F)-(G_{n+1})_{[x]}$, or, equivalently, the images of $g.x$ and $x$ in $\sB(G_{n+1}^{\rm{der}}, k)$ are distinct. 
Let 
$f$ be an element of $\Hom_{\KYu \cap {^{g}\KYu}}\left(^{g}\repKYu, {\repKYu}\right) - \{0\}$.  
We denote its image in the space $V_{\repKYu}$ of the representation of $\repKYu$ by $V_f$. 
We write $H_{n+1}$ for the derived subgroup $G_{n+1}^{\text{der}}$ of $G_{n+1}$ and denote by $(H_{n+1})_{x,r}$ the Moy--Prasad filtration subgroup of depth $r \in \bR_{\geq 0}$ at the image of $x$ in $\sB(H_{n+1},F)$. Then 
$^{g}(H_{n+1})_{x,0}=(H_{n+1})_{g.x,0}$, and we have 
\begin{equation} \label{eqn-intertwiner3} f \in \Hom_{\KYu \cap {^{g}\KYu}}\left(^{g}\repKYu, {\repKYu}\right) - \{0\} \subset \Hom_{(H_{n+1})_{x,0} \cap {(H_{n+1})_{g.x,0}}}\left(^{g}\repKYu, {\repKYu}\right) - \{0\} \, . 
\end{equation}

By construction  $\repKYu|_{(H_{n+1})_{x,0+}}=\hat \phi|_{(H_{n+1})_{x,0+}} \cdot \Id = \phi|_{(H_{n+1})_{x,0+}} \cdot \Id$ where  $ \hat \phi := \prod_{i=1}^n \hat \phi_i|_{(G_{n+1})_{[x]}G_{x,\frac{r_1}{2}+}}$ and  $\phi:=\prod_{i=1}^n \phi_i|_{G_{n+1}(F)}$. In addition, for all $y \in (H_{n+1})_{g.x,0+}$ we have
$^{g}\repKYu(y)=\hat\phi(g^{-1}yg)=\prod_{i=1}^n \phi_i(g^{-1}yg)=\prod_{i=1}^n \phi_i(y)=\phi(y) $, because $g \in G_{n+1}(F)$. 
Hence, by \eqref{eqn-intertwiner3}, the action of 
$$U:=((H_{n+1})_{x,0} \cap {(H_{n+1})_{g.x,0+}})(H_{n+1})_{x,0+}$$
  on the image $V_f$ of $f$ via ${\repKYu}$ is given by $\phi \cdot \Id$.

Recall that the image of $x$ in $\sB(H_{n+1}, k)$ is a vertex by Condition \eqref{condition-vertex} of the input in Section~\ref{Section-input}. Hence the group $(((H_{n+1})_{x,0} \cap {(H_{n+1})_{g.x,0+}})(H_{n+1})_{x,0+})/(H_{n+1})_{x,0+}$ is the ($\bF_q$-points of) a unipotent radical of a (proper) parabolic subgroup of $(H_{n+1})_{x,0}/(H_{n+1})_{x,0+}$. We denote this subgroup by $\Uff$.

In the remainder of the proof we exhibit a subspace $V_\kappa' \subset V_\kappa$ such that $V_f \subset  V_\rho \otimes V_\kappa'$ and prove that the action of $\UU$ on $V_\kappa'$  via $\kappa$ is given by $\phi \cdot \Id$. Hence, since $\UU$ also acts via $\phi \cdot \Id$ on $V_f \subset  V_\rho \otimes V_\kappa'$, we deduce that $(\rho|_{\Uff}, V_\rho)$ contains the trivial representation, which contradicts that $\rho|_{(G_{n+1})_{x,0}}$ is cuspidal (see Condition \eqref{condition-cuspidal} of the input in Section \ref{Section-input}).

Let $T$ be a maximal torus of $G$ that splits over a tamely ramified extension $E$ of $F$ such that $x$ and $g.x$ are contained in $\sA(T,E)$. (Such a torus exists by Remark \ref{Rem-tame-torus}.) Let $\lambda \in X_*(T) \otimes_\bZ\bR=\Hom_{\ov F}(\bG_m,T_{\ov F}) \otimes_\bZ \bR$ such that $g.x = x + \lambda$, and observe that $\Uff$ is the image of 
\begin{equation}
	U_{n+1}:=H_{n+1}(F) \cap \<U_\alpha(E)_{x,0} \, | \, \alpha \in \Phi(G_{n+1}, T), \lambda(\alpha) >0 \>
\end{equation}
in $(H_{n+1})_{x,0}/(H_{n+1})_{x,0+}$.
 We define for $1 \leq i \leq n$
$$U_i:= G(F) \cap \<U_\alpha(E)_{x,\frac{r_i}{2}} \, | \, \alpha \in \Phi(G_i, T)-\Phi(G_{i+1}, T), \lambda(\alpha) >0 \> .$$

Note that ${^{g}(G_i)_{x,r_i,\frac{r_i}{2}+}}=(G_i)_{g.x,r_i,\frac{r_i}{2}+}$ and $U_i \subset (G_i)_{x,r_i,\frac{r_i}{2}} \cap (G_i)_{g.x,r_i,\frac{r_i}{2}+} \subset \KYu \cap {^{g}\KYu}$. More precisely,  $U_i \subset (G_i)_{g.x,r_i+,\frac{r_i}{2}+}$, hence ${^{{g}^{-1}}U_i} \subset (G_i)_{x,r_i+,\frac{r_i}{2}+}$ and $\hat \phi_j|_{{^{{g}^{-1}}U_i}}$ is trivial for $j \geq i$. Thus, by Equation \eqref{eqn-g'-rho}, we obtain that 
$^{g}\repKYu|_{U_i} = \prod_{j=1}^{i-1} \hat\phi_j \cdot \Id$. Hence $U_i$ acts on $V_f$ via the character $\prod_{j=1}^{i-1} {\hat\phi_j}|_{U_i}=\prod_{1 \leq j \leq n \atop j \neq i} {\hat\phi_j}|_{U_i}$. Since $U_i$ acts trivially via $\rho$ on the space $V_\rho$ underlying the representation of $\rho$ and $U_i$ acts via $\prod_{1 \leq j \leq n \atop j \neq i} {\hat\phi_j}|_{U_i}$ on $\bigotimes_{1 \leq j \leq n \atop j \neq i} V_{\omega_j}$, we obtain
 $$V_f \subset V_\rho \otimes \bigotimes_{j=1}^{i-1} V_{\omega_j} \otimes V_{\omega_i}^{U_i} \otimes \bigotimes_{j=i+1}^{n} V_{\omega_j} \subset V_\rho \otimes \bigotimes_{i=1}^n V_{\omega_i} = V_\rho \otimes V_\kappa$$ 
 for $1 \leq i \leq n$. Hence we deduce that $V_f \subset V_\rho \otimes \bigotimes_{i=1}^n V_{\omega_i}^{U_i}$. We will see that $\bigotimes_{i=1}^n V_{\omega_i}^{U_i}$ is the subspace $V_\kappa' \subset V_\kappa$ that we are looking for.

In order to study the subspace $V_{\omega_i}^{U_i}$ for $1 \leq i \leq n$, we recall that we write $V_i=(G_i)_{x,r_i,\frac{r_i}{2}}/(G_i)_{x,r_i,\frac{r_i}{2}+}$ and equip $V_i$ with the pairing $\<\cdot, \cdot\>_i$ defined by $\<a,b\>_i=\hat \phi_i(aba^{-1}b^{-1})$ (using the above fixed identification of the $p$-th roots of unity in $\bC^*$ with $\bF_p$). We define the space
$V_i^+$ to be the image 
of $U_i= G(F) \cap \<U_\alpha(E)_{x,\frac{r_i}{2}} \, | \, \alpha \in \Phi(G_i, T)-\Phi(G_{i+1}, T), \lambda(\alpha)>0 \> $ in $V_i$, 
the space
$V_i^0$ to be the image 
of $ G(F) \cap \<U_\alpha(E)_{x,\frac{r_i}{2}} \, | \, \alpha \in \Phi(G_i, T)-\Phi(G_{i+1}, T), \lambda(\alpha)=0 \> $ in $V_i$, 
and 
$V_i^-$ to be the image 
of $ G(F) \cap \<U_\alpha(E)_{x,\frac{r_i}{2}} \, | \, \alpha \in \Phi(G_i, T)-\Phi(G_{i+1}, T), \lambda(\alpha)<0 \> $ in $V_i$.
Then $V_i=V_i^+ \oplus V_i^{0} \oplus V_i^-$ and the $\bF_p$-vector subspaces $V_i^+$ and $V_i^-$ are both totally isotropic. Since $\phi_i$ is $G_i$-generic of depth $r_i$ relative to $x$ the orthogonal complement of $V_i^+$ is $V_i^+ \oplus V_i^0$, the orthogonal complement of $V_i^-$ is $V_i^0 \oplus V_i^-$, and $V_i^0$ is a non-degenerate subspace of $V_i$. 
We denote by $P_i \subset \Sp(V_i)$ the (maximal) parabolic subgroup of $\Sp(V_i)$ that preserves the subspace $V_i^+$ and that therefore also preserves $V_i^+ \oplus V_i^0$. We obtain a surjection $\pr_{i,0}: P_i \twoheadrightarrow \Sp(V_i^0)$ by composing restriction to $V_i^0$ with projection from $V_i^+ \oplus V_i^0$ to $V_i^0$ with kernel $V_i^+$.  Note that the image $\Uff_{\Sp(V_i)}$ of $\Uff$ in $\Sp(V_i)$ is contained in $P_i$ and that $\pr_{i,0}(\Uff_{\Sp(V_i)})=\Id_{V_i^0}$. 

Recall that $V_i^\sharp$ is the Heisenberg group with underlying set $V_i \times \bF_p$ that is attached to the symplectic $\bF_p$-vector space  $V_i$ with pairing $\< \cdot, \cdot \>_i$, 
 and note that the subset $V_i^0 \times \bF_p \subset V_i \times \bF_p$ forms a subgroup, which is the Heisenberg group $(V_i^0)^\sharp$ attached to the symplectic vector space  $V_i^0$ with the (restriction of the) pairing $\< \cdot, \cdot \>_i$. We denote by $V_{\omega_i}^0$ a Weil--Heisenberg representation of $\Sp(V_i^0) \ltimes (V_i^0)^\sharp$ corresponding to the same central character as the central character of $V_i^\sharp$ acting on $V_{\omega_i}$ (which in turn corresponds to the character $\hat \phi_i|_{(G_i)_{x,r_i,\frac{r_i}{2}+}}$ via the special isomorphism $j_i$).
By \cite[Theorem~2.4.(b)]{Gerardin}  
the restriction of the Weil--Heisenberg representation $V_{\omega_i}$ from $\Sp(V_i) \ltimes V_i^\sharp$ to $P_i \ltimes V_i^\sharp$ is 
given by 
$$\Ind_{P_i \ltimes ( V_i^+ \times (V_i^0)^\sharp)}^{P_i \ltimes V_i^\sharp} V_{\omega_i}^0 \otimes (\bC_{\chi^{V_i^+}} \ltimes 1),$$
where the group $P_i \ltimes ( V_i^+ \times (V_i^0)^\sharp)$ acts on $V_{\omega_i}^0$ by composing the projection 
$$\pr_{i,0} \ltimes (\pr_{+0,0} ): P_i \ltimes ( V_i^+ \times (V_i^0)^\sharp) \ra \Sp(V_i^0) \ltimes (V_i^0)^\sharp$$
 (where $\pr_{+0,0}: ( V_i^+ \times (V_i^0)^\sharp) \twoheadrightarrow (V_i^0)^\sharp$ denotes the projection with kernel $V_i^+$) with the Weil--Heisenberg representation of $\Sp(V_i^0) \ltimes (V_i^0)^\sharp$, and $\bC_{\chi^{V_i^+}}$ is a one dimensional space on which the action of $P_i$ is given by a quadratic character\footnote{The definition of $\chi^{V_i^+}$ is $\det(\pr_{i,+}(\cdot))^{(p-1)/2}$, but we will not need the precise definition for our proof. Note that the statement of \cite[Theorem~2.4.(b)]{Gerardin} omits the factor $\chi^{V_i^+} \ltimes 1$ (denoted by $\chi^{E_+}$ in \cite[Theorem~2.4.(b)]{Gerardin}), which is a typo that was pointed out by Loren Spice.} $\chi^{V_i^+}$ that factors through the projection $\pr_{i,+}: P_i \rightarrow \GL(V_i^+)$ obtained by restricting elements in $P_i$ to $V_i^+$. 

Let $\Uiff$ be the image of $U_i$ in the Heisenberg group $(G_i)_{x,r_i,\frac{r_i}{2}}/\left((G_i)_{x,r_i,\frac{r_i}{2}+} \cap \ker(\hat \phi_i)\right)$. Then by Yu's construction of the special isomorphism $j_i:(G_i)_{x,r_i,\frac{r_i}{2}}/\left((G_i)_{x,r_i,\frac{r_i}{2}+} \cap \ker(\hat \phi_i)\right) \ra V_i^\sharp$
in \cite[Proposition~11.4]{Yu},
we have $j_i(\Uiff)=V_i^+ \times 0 \subset V_i \times \bF_p$. 
Since the orthogonal complement of $V_i^-$ is $V_i^0 \oplus V_i^-$, and hence for every element $v_- \in V_i^-$ there exists $v_+ \in V_i^+$ such that $\<v_-,v_+\>_i \neq 0$, we have
\begin{equation} \label{eqn-fixed-points}\left(\Ind_{P_i \ltimes ( V_i^+ \times (V_i^0)^\sharp)}^{P_i \ltimes V_i^\sharp} V_{\omega_i}^0 \otimes (\bC_{\chi^{V_i^+}} \ltimes 1)\right)^{1 \ltimes (V_i^+ \times 0)} \simeq V_{\omega_i}^{0}  \otimes (\bC_{\chi^{V_i^+}} \ltimes 1)
\end{equation}
 as a representation of $P_i$. 

Note that the image of $\Uff_{\Sp(V_i)}$ in $\GL(V_i^+)$ under the projection $\pr_{i,+}:P_i \ra \GL(V_i^+)$ is unipotent since $\Uff$ is unipotent. Hence $\pr_{i,+}(\Uff_{\Sp(V_i)})$ is contained in the commutator subgroup of $\GL(V_i^+)$, and $\chi^{V_i^+}|_{\pr_{i,+}(\Uff_{\Sp(V_i)})}$ is trivial. Moreover, we observed above that $\pr_{i,0}(\Uff_{\Sp(V_i)})=\Id_{V_i^0}$. Thus $\Uff_{\Sp(V_i)}$ acts trivially on $V_{\omega_i}^{0}  \otimes (\bC_{\chi^{V_i^+}} \ltimes 1)$. 

Recall that the action of $\UU$ on $(V_{\omega_i})^{U_i}$ is given by the product of $\phi_i|_\UU$ with the above Weil representation construction, see Section \ref{Section-construction}, page \pageref{page-tildephi}. Hence $\UU$ acts on $(V_{\omega_i})^{U_i}$ via the character $\phi_i|_\UU$. Since we proved above that $\UU$ acts via $\phi=\prod_{i=1}^n \phi_i$ on $V_f \subset V_\rho \otimes \bigotimes_{i=1}^n (V_{\omega_i})^{U_i}$, we deduce that there exists a non-trivial subspace $V_{\rho,f}$ of $V_\rho$ on which $\UU$ acts trivially. Hence $\rho|_{\Uff}$ contains the trivial representation, which contradicts that $\rho|_{(G_{n+1})_{x,0}}$ is cuspidal. \qed

\section{A counterexample} \label{Section-counterexample}
In this section we provide a counterexample to \cite[Proposition~14.1~and~Theorem~14.2]{Yu}, whose proof relied on the misprinted version of \cite[Theorem~2.4(b)]{Gerardin}. To state the content of the section more precisely, let $G'$ be a tamely ramified twisted Levi subgroup of $G$, let $x \in \sB(G', F)$, and $\phi$ a character of $G'(F)$ that is $G$-generic relative to $x$ of depth $r$ for some $r \in \bR_{>0}$, i.e. we are in the setting of \cite[\S~14]{Yu}. Following \cite{Yu}, we set 
\begin{eqnarray*}
&& J=(G', G)(F)_{x,(r,\frac{r}{2})}, \,  J_+=(G', G)(F)_{x,(r,\frac{r}{2}+)}, \\
&& K=G'(F) \cap G_{[x]}, \,  
 K_+=G'(F) \cap G_{x,0+}, \,  N=\ker \hat \phi 
 \end{eqnarray*} 
where $\hat \phi$ is defined as in \cite[\S4 and \S9]{Yu}, see also page \pageref{page-hatphi} of this paper, and we denote by $\tildephi$ the representation of $K \ltimes J$ which is the pull back of the Weil representation of $\Sp(J/J_+) \ltimes (J / N)$ via the symplectic action given by \cite[Proposition~11.4]{Yu}, see also page \pageref{page-tildephi} of this paper.

In this section, we provide an example for $G' \subset G$, $x$ and $\phi$ as above and $g \in G'(F)$ such that 
\begin{equation} \label{eqn-Yu-counterexample}
\dim\Hom_{(K \cap {^{g}K}) \ltimes {(J \cap {^{g}J})}}({^g\tildephi}, \tildephi) =0.
\end{equation}
Following \cite[\S14]{Yu} we denote by $\phi'$ the representation of $KJ$ whose inflation $\inf \phi'$ to $K \ltimes J$ yields $\inf(\phi|_K) \otimes \wt \phi$.
By the discussion in \cite{Yu} immediately following Theorem 14.2 (see also Corollary \ref{Cor-Yu-counterexample} below), Equation \eqref{eqn-Yu-counterexample} implies that $$\dim \Hom_{KJ \cap {^g(KJ)}}(^g\phi',\phi')=0$$ and therefore provides a counterexample to the claim that $\Hom_{KJ \cap {^g(KJ)}}(^g\phi',\phi')$ has always dimension one that was made in \cite[Proposition~14.1]{Yu} and in its more general version \cite[Theorem~14.2]{Yu}.

Consider the case $G=\Sp_{10}$ over $F$ corresponding to the symplectic pairing given by 
$ \begin{pmatrix} 
0 & J_5 \\  
-J_5 & 0 \\
\end{pmatrix}$ where
$$ J_5 = \begin{pmatrix} 
0 & 0 & 0 & 0 & 1 \\  
0 & 0 & 0 & -1 & 0 \\
0 & 0 & 1 & 0 & 0 \\
0 & -1 & 0 & 0 & 0 \\
1 & 0 & 0 & 0 & 0 \\
\end{pmatrix} .$$
We assume that the residue field of $F$ is $\bF_p$ for some prime number  $p>10$.

Let $T \subset \Sp_{10}$ be the diagonal maximal torus using the standard coordinates, and write $\ft=\Lie(T)(F)$.
We identify the apartment $\sA(T, F)$ with $X_*(T)\otimes \bR$ ($X_*$ being as above the cocharacters over $\ov F$, or, equivalently, the cocharacters over $F$) using the standard parametrization of the root groups as base point, i.e. the point for which the attached parahoric subgroup is $\Sp_{10}(\cO)$ in the standard coordinates.
 Identifying $X_*(T) \otimes_\bZ \bR$ with $\bR^5$ where the first standard basis vector corresponds to $t \ra \diag(t, 1, 1, 1, 1, 1, 1, 1, 1, t^{-1})$, the second to $t \ra \diag(1, t, 1, 1, 1, 1, 1, 1, t^{-1}, 1)$, etc., we let $x$ be the point of $\sA(T, F)$ corresponding to $(-\frac{1}{4}, 0, 0, \frac{1}{4}, \frac{1}{4})$.

Let $\pi$ be a uniformizer of $F$,  and let $\varpi$ in $F^{sep}$ such that $\varpi^2=\pi$.
 Let $X \in \fg^*_{x,-\frac{1}{2}}$ be the element given by
 $$ \fsp_{10} \ni (A_{i,j}) \mapsto \pi^{-1}A_{1,10}+A_{10,1} .$$
  We set $G'$ to be the centralizer $\Cent_{G}(X)$ of $X$ in $G$. Note that 
$$G'=\Cent_G\left(\begin{pmatrix} 
0 & 0_{1 \times 8} & 1 \\  
0_{8 \times 1} & 0_{8\times 8} & 0_{8 \times 1} \\
\pi^{-1} & 0_{1 \times 8} & 0 \\  
\end{pmatrix} \right) \simeq \U(1) \times \Sp_8$$
 is a twisted Levi subgroup of $G=\Sp_{10}$ (with anisotropic center).  
\begin{Lemma}
	The (restriction to $\fg'$ of the) element $X$ is $G$-generic of depth $r=\frac{1}{2}$ (for the pair $G'\subset G$), and the point $x=(-\frac{1}{4}, 0, 0, \frac{1}{4}, \frac{1}{4}) \in \sA(T,F)\subset \sB(G,F)$ is contained in $\sB(G',F)$.
\end{Lemma}
\Proof
First note that $X$ is $G'$-invariant by construction.
Let  $\sqrt{2\varpi}$ in $F^{sep}$ such that $\sqrt{2\varpi}^2=2 \varpi$. We consider the maximal torus 
$$T'= 
\begin{pmatrix} 
 \frac{\varpi}{\sqrt{2\varpi}} & 0_{1 \times 8} & \frac{-\varpi }{\sqrt{2\varpi}}\\  
0_{8 \times 1} & 1_{8\times 8} & 0_{8 \times 1} \\
\frac{1}{\sqrt{2\varpi}} & 0_{1 \times 8} & \frac{1}{\sqrt{2\varpi}} \\  
\end{pmatrix}
T
\begin{pmatrix} 
\frac{1}{{\sqrt{2\varpi}}} & 0_{1 \times 8} & \frac{\varpi}{\sqrt{2\varpi}} \\  
0_{8 \times 1} & 1_{8\times 8} & 0_{8 \times 1} \\
\frac{-1}{{\sqrt{2\varpi}}} & 0_{1 \times 8} & \frac{\varpi}{\sqrt{2\varpi}} \\  
\end{pmatrix}
 $$
 of $G'$. 
Then the set $\{H_{\alpha}=d\check\alpha(1) \, | \, \alpha \in \Phi(G, T') \setminus \Phi(G',T')\}$, where $\check \alpha$ denotes the dual root of $\alpha$ 
is given by
$$ \pm \begin{pmatrix} 
0 & 0_{1 \times 8} & \varpi \\  
0_{8 \times 1} & 0_{8\times 8} & 0_{8 \times 1} \\
\varpi^{-1} & 0_{1 \times 8} & 0 \\  
\end{pmatrix} $$
and sums of the former with a diagonal matrix of the form
\begin{eqnarray*}
\pm \diag(0,1,0,0,0,0,0, 0,-1,0), \,
\pm \diag(0,0,1,0,0,0,0,-1,0,0), \\
\pm \diag(0,0,0,1,0,0,-1,0,0,0), \,
\pm \diag(0,0,0,0,1,-1,0,0,0,0) \, .
\end{eqnarray*}
Hence $\val(X(H_\alpha))=\val(\pm2\varpi^{-1})=-\frac{1}{2}$ for all $\alpha \in \Phi(G, T') \setminus \Phi(G',T')$, where $\val$ denotes the valuation of $F$ with image $\bZ\cup\{\infty\}$.
Moreover, $X$ is contained in $(\fg')^*_{y,-\frac{1}{2}}$ for any point $y \in \sB(G', F)$ and its  restriction to the Lie algebra of the identity component $Z(G')^\circ$ of the center of $G'$ is contained in $\Lie^*(Z(G')^\circ)_{-\frac{1}{2}}$. Since $p>10$, i.e. $p$ does not divide the order of the Weyl group of $\Sp_{10}$, we conclude that $X$ is $G$-generic of depth $\frac{1}{2}$. 

Let $E$ be a finite, tamely ramified extension of $F$ that contains $\sqrt{2\varpi}$. Then $x$ is a vertex in $\sB(G,E)$ and it is an easy calculation to check that 
$$G_{x,0}(E)=\begin{pmatrix} 
\frac{\varpi}{\sqrt{2\varpi}} & 0_{1 \times 8} & \frac{-\varpi }{\sqrt{2\varpi}}\\  
0_{8 \times 1} & 1_{8\times 8} & 0_{8 \times 1} \\
\frac{1}{\sqrt{2\varpi}} & 0_{1 \times 8} & \frac{1}{\sqrt{2\varpi}} \\  
\end{pmatrix}
G_{y,0}(E)
\begin{pmatrix} 
\frac{1}{{\sqrt{2\varpi}}} & 0_{1 \times 8} & \frac{\varpi}{\sqrt{2\varpi}} \\  
0_{8 \times 1} & 1_{8\times 8} & 0_{8 \times 1} \\
\frac{-1}{{\sqrt{2\varpi}}} & 0_{1 \times 8} & \frac{\varpi}{\sqrt{2\varpi}} \\  
\end{pmatrix}$$
for the point $y=(0,0,0,\frac{1}{4},\frac{1}{4}) \in \sA(T,F) \subset \sA(T,E)$. Hence $x=\begin{pmatrix} 
\frac{\varpi}{\sqrt{2\varpi}} & 0_{1 \times 8} & \frac{-\varpi }{\sqrt{2\varpi}}\\  
0_{8 \times 1} & 1_{8\times 8} & 0_{8 \times 1} \\
\frac{1}{\sqrt{2\varpi}} & 0_{1 \times 8} & \frac{1}{\sqrt{2\varpi}} \\  
\end{pmatrix}.y$, which implies $x \in \sA(T', E) \subset \sB(G', E)$, and therefore $x \in \sB(G', E) \cap \sB(G, F) =\sB(G',F)$.
 \qed

The element $X$ yields a linear map from $\fg'_{x,\frac{1}{2}}$ to $\cO$ that sends $\fg'_{x,\frac{1}{2}+}$ to $\varpi\cO$ and defines a character of $G'_{x,\frac{1}{2}}$ that is trivial on $G'_{x,\frac{1}{2}+}$ and trivial on $G'_{x,\frac{1}{2}} \cap \Sp_8 \subset G'_{x,\frac{1}{2}} \cap (U(1) \times \Sp_8)(F) \simeq G'_{x,\frac{1}{2}} $. Since $U(1)$ is abelian, we can extend this character to a character of $G'(F)$ (trivial on $\Sp_8(F) \subset G'(F)$), which we denote by $\phi$. Since $X$ is $G$-generic of depth $r=\frac{1}{2}$ (for the pair $G'\subset G$), the character $\phi$ is $G$-generic relative to $x$ of depth $r$ in the sense of Yu (\cite[\S~9]{Yu}).

\begin{Prop}\label{Prop-Yu-counterexample}
	Let $G'\subset G, \phi, x, r$ as above and 
		 $$g= \begin{pmatrix} 
		1&  0 & 0 & 0 \\  
		0&  J_4^+ & 0 &0 \\  
		0 & 0 & -J_4 ^+ &0 \\
		0&  0 & 0 & 1 \\  
		\end{pmatrix} \in G'(F) \text{ where }
		 J_4^+ = \begin{pmatrix} 
		0 & 0 & 0 & 1 \\  
		0 & 0 & 1 & 0 \\
		0 & 1 & 0 & 0 \\
		1 & 0 & 0 & 0 \\
		\end{pmatrix} . $$
		Then 
	$\dim\Hom_{(K \cap {^{g}K}) \ltimes {(J \cap {^{g}J})}}({^g\tildephi}, \tildephi) =0. $
\end{Prop}
\Proof
Using the standard coordinates we define the groups
$$ \Htt:=   \begin{pmatrix} 
1 & 0 & 0 & 0 & 0 \\  
0 & \GL_2(\cO) & 0 & 0 & 0  \\
0 & 0 & 1_{4 \times 4} & 0 & 0  \\
0 & 0 & 0 & \GL_2(\cO) & 0 \\
0 & 0 & 0 & 0 & 1  \\
\end{pmatrix} \cap \Sp_{10}(F) \subset G'_{x,0}$$ 
$$ \Hff:=   \begin{pmatrix} 
1_{3\times3} & 0 & 0 & 0 \\  
0 & \GL_2(\cO) & 0 & 0  \\
0 & 0 & \GL_2(\cO) & 0 \\
0 & 0 & 0 & 1_{3\times3} \\
\end{pmatrix} \cap \Sp_{10}(F) \subset G'_{x,0}$$

Note that $\Htt \simeq \GL_2(\cO) \simeq \Hff$, and $g\Htt g^{-1}=\Hff$ and $g\Hff g^{-1}=\Htt$, hence $\Htt \in K \cap {^gK}$. Moreover, the image of $\Htt$ and the image of $\Hff$ in $G'_{x,0}/G'_{x,0+}$ are both isomorphic to $\GL_2(\bF_p)$.

We are going to show that $\dim\Hom_{\Htt}({^g\tildephi}, \tildephi)=0$, which implies that \newline
$\dim\Hom_{({K\cap{^{g}K}}) \ltimes ({J \cap {^{g}J}})}({^g\tildephi}, \tildephi) = 0$. 

We write $V=J/J_+$ where we recall that $J=(G',G)(F)_{x,(\frac{1}{2},\frac{1}{4})}$ and  $J_+=(G',G)(F)_{x,(\frac{1}{2},\frac{1}{4}+)}$.
Let $\fg(\cO)$ be the $\cO$-points of the Lie algebra of the reductive parahoric group scheme over $\cO$ corresponding to the base point $(0,0,0,0,0)$, i.e. the Lie algebra of $\Sp_{10}$ defined over $\cO$ in the standard basis. We denote by $\fg(\cO)_{t_1t_2^{-1}}$ the submodule of $\fg(\cO)$ corresponding to the root $\diag(t_1, t_2, t_3, t_4, t_5, t_5^{-1}, t_4^{-1}, t_3^{-1}, t_2^{-1}, t_1^{-1}) \mapsto t_1t_2^{-1}$ and analogously for all other indices.
Then the four-dimensional $\bF_p$-vector space $V$ is spanned by the images of 
\begin{equation*} \fg(\cO)_{t_1^{-1}t_2} , \quad
	 \fg(\cO)_{t_1^{-1}t_3} , \quad
	 \fg(\cO)_{t_1^{-1}t_2^{-1}} , \quad
	 \fg(\cO)_{t_1^{-1}t_3^{-1}} .
\end{equation*} 
Each of these images is a one dimensional $\bF_p$-vector subspace of $V$, which we denote by $V_{t_2}$, $V_{t_3}$, $V_{t_2^{-1}}$ and $V_{t_3^{-1}}$, respectively.
The pairing on $V=J/J_+$ defined by $\<a,b\>=\hat \phi(aba^{-1}b^{-1})$ for $a, b \in J$ turns $V$ into a symplectic $\bF_p$-vector space, and $V^+:=V_{t_2} \oplus V_{t_3}$ and  $V^-:=V_{t_2^{-1}} \oplus V_{t_3^{-1}}$ are both maximal isotropic subspaces.
Recall that $\tildephi$ is defined to be the pullback  to $K\ltimes J$  of the Weil--Heisenberg representation of $\Sp(V) \ltimes (J/N)$ via the symplectic action defined in \cite[Proposition~11.4]{Yu}.
Hence the actions of $\Htt$ and $\Hff$ on $\tildephi$ factor through $\Sp(V)$, and therefore $\Hff$ acts trivially on $\tildephi$. Thus $^g\tildephi|_{\Htt}={^g\tildephi}|_{g\Hff g^{-1}}$ is trivial, and in order to prove that $\dim\Hom_{\Htt}({^g\tildephi}, \tildephi)=0$ it suffices to show that the representation $\tildephi|_{\Htt}$ has no non-zero $\Htt$-fixed vector.

We denote by $P$ the parabolic subgroup of $\Sp(V)$ that preserves $V^+$. Then the image of $\Htt$ in $\Sp(V)$ is the Levi subgroup $M\simeq \GL(V^+)\simeq \GL(V^-)$ of $P$ that stabilizes $V^+$ and $V^-$. Recall that we denote by $V^\sharp$ the group with underlying set $V \times \bF_p$ and with group law $(v,a).(v',a')=(v+v', a+a'+\frac{1}{2}\<v,v'\>)$.
By \cite[Theorem~2.4.(b)]{Gerardin}  
the restriction of the Weil--Heisenberg representation from $\Sp(V) \ltimes V^\sharp$ to $P \ltimes V^\sharp$ is 
given by 
$$\pi:= \Ind_{P \ltimes (V^+ \times \bF_p)}^{P \ltimes V^\sharp} \chi^{V^+} \ltimes \phi,$$
where $\chi^{V^+}$ is the character\footnote{Note that the statement of \cite[Theorem~2.4.(b)]{Gerardin} omits the character $\chi^{V^+}$ in the induction, which is a typo that was pointed out by Loren Spice.} of $P$ given by $P \ni p \mapsto \det(p|_{V^+})^{\frac{p-1}{2}} \in \{ \pm 1\} \subset \bC^*$ and we denote by abuse of notation by $\phi$ the (restriction to $V^+ \times \bF_p$ of the) character $\phi \circ j^{-1}$, where $j:J/N \xrightarrow{\simeq} V^\sharp$ denotes the special isomorphism from \cite[Proposition~11.4]{Yu}.

Let $f: P \ltimes V^\sharp \ra \bC$ be an element of the representation space of $\pi$ 
and suppose that $f $ is non-zero and $M$-invariant. Hence there exists $v \in V^-$ such that $f(1 \ltimes v) \neq 0$. Let $v' \in V^-$ so that $v$ and $v'$ form a basis of $V^-$ and let $m=\begin{pmatrix} 1 & 0 \\ 0 & a \end{pmatrix} \in \GL(V^-)$ using the basis $(v, v')$ where $a \in \bF_p$ such that $a^{\frac{p-1}{2}}=-1$. Identifying $\GL(V^-)$ with $M$ (via the action of $M$ on $V^-$), we obtain that 
\begin{eqnarray*} 
	f(1\ltimes v) &=& m.f(1 \ltimes v)=f((1 \ltimes v)(m \ltimes 1))=f((m \ltimes 1)(1 \ltimes m^{-1}.v)) \\
	&=&\chi^{V^+}(m)f(1 \ltimes m^{-1}.v) 
	= \det(m|_{V^+})^{\frac{p-1}{2}}f(1 \ltimes v)=\det(m|_{V^-})^{\frac{p-1}{2}}f(1 \ltimes v)\\
	&=&-f(1 \ltimes v)
\end{eqnarray*}
This contradicts that $f(1 \ltimes v) \neq 0$, hence the representation $\pi$ does not contain any non-zero element fixed under the action of $M$.
Therefore $\tildephi$ does not contain any non-zero element fixed under the action of $\Htt$. 
Thus $\dim\Hom_{\Htt}({^g\tildephi}, \tildephi)=0$. 
\qed

\begin{Cor} \label{Cor-Yu-counterexample}
	In the setting of Proposition \ref{Prop-Yu-counterexample}, we have
	$$\dim \Hom_{KJ \cap {^g(KJ)}}(^g\phi',\phi')=0 . $$
\end{Cor}
\Proof This follows from the fact that $KJ \cap {^g(KJ)}=(K \cap {^{g}K})(J \cap {^{g}J})$ (\cite[Lemma~13.7]{Yu}) as discussed in \cite{Yu} in the lines  immediately following Theorem 14.2.
\qed

\bibliography{Fintzenbib}  

\end{document}